\documentclass[11pt,a4paper]{article}

\usepackage{amsmath,amsthm,amstext,amscd,amssymb,euscript,mathrsfs}
\usepackage{calrsfs}
\usepackage{epsfig}
\usepackage{color}

\newcommand{\Z}{\mathbb Z}
\newcommand{\R}{\mathbb R}
\newcommand{\N}{\mathbb N}

\newcommand{\E}{\mathbb E}
\newcommand{\Zd}{\mathbb Z^d}

\newcommand{\epsi}{\ensuremath{\epsilon}}

\newcommand{\pee}{\ensuremath{\mathbb{P}}}

\newcommand{\loc}{\mathcal{L}}

\newcommand{\x}{{\bf x}}
\def\1{{\mathchoice {\rm 1\mskip-4mu l} {\rm 1\mskip-4mu l}
{\rm 1\mskip-4.5mu l} {\rm 1\mskip-5mu l}}}

\newtheorem{theorem}{{\small T}{\scriptsize HEOREM}}[section]
\newtheorem{corollary}{{\bf{\small C}{\scriptsize OROLLARY}}}[section]
\newtheorem{proposition}{{\bf{\small P}{\scriptsize ROPOSITION}}}[section]
\newtheorem{lemma}{{\bf{\small L}{\scriptsize EMMA}}}[section]
\newtheorem{remark}{{\bf{\small R}{\scriptsize EMARK}}}[section]
\newtheorem{definition}{{\bf{\small D}{\scriptsize EFINITION}}}[section]

\renewenvironment{proof}[1]
{\noindent{{\bf{\small{ P}{\scriptsize ROOF}}}.}\hspace{0.1cm} #1} {$\;\qed$\newline}

\newcommand{\beq}{\begin{eqnarray}}
\newcommand{\eeq}{\end{eqnarray}}

\newcommand{\ba}{\begin{align*}}
\newcommand{\ea}{\end{align*}}

\newcommand{\be}{\begin{equation}}
\newcommand{\ee}{\end{equation}}

\newcommand{\bl}{\begin{lemma}}
\newcommand{\el}{\end{lemma}}

\newcommand{\br}{\begin{remark}}
\newcommand{\er}{\end{remark}}

\newcommand{\bt}{\begin{theorem}}
\newcommand{\et}{\end{theorem}}

\newcommand{\bd}{\begin{definition}}
\newcommand{\ed}{\end{definition}}

\newcommand{\bp}{\begin{proposition}}
\newcommand{\ep}{\end{proposition}}

\newcommand{\bc}{\begin{corollary}}
\newcommand{\ec}{\end{corollary}}

\newcommand{\bpr}{\begin{proof}}
\newcommand{\epr}{\end{proof}}

\newcommand{\bi}{\begin{itemize}}
\newcommand{\ei}{\end{itemize}}

\newcommand{\ben}{\begin{enumerate}}
\newcommand{\een}{\end{enumerate}}

\newcommand{\caB}{{\mathcal B}}

\newcommand{\caD}{{\EuScript D}}

\newcommand{\caL}{{\mathcal L}}

\newcommand{\caV}{{\mathcal V}}

\newcommand{\mee}{\frac{m}{2}}

\renewcommand{\L}{\mathcal{L}}

\begin{document}
\title{Coupling and hydrodynamic limit for the inclusion process}
\author{ Alex Opoku and
Frank Redig\\
{\small Delft Institute of Applied Mathematics,}\\
{\small Technische Universiteit Delft}\\
{\small Mekelweg 4, 2628 CD Delft, Nederland}\\
}

\maketitle

\begin{abstract}
We show propagation of local equilibrium for the symmetric inclusion process (SIP) after diffusive rescaling of space and time, as well as the local equilibrium property of the non-equilibrium steady state in the boundary driven SIP. The main tool is self-duality and a coupling between
$n$ SIP particles and $n$ independent random walkers.
\bigskip

\noindent

\end{abstract}

\section{Introduction}
The symmetric inclusion process (SIP) was introduced in \cite{GKRV} and studied in more detail in \cite{GRV}.
It is a process where particles move on the lattice $\Zd$ according to symmetric nearest neighbor random walk jumps at rate
$m/2$, and on top of that  there is an attractive interaction governed by ``inclusion jumps'' where particles at nearest neighbor positions
invite each other (invitations which are always followed up) at rate 1. The attractive interaction as opposed to exclusion in the well-known
exclusion process is responsible for interesting phenomena. In the limit of small diffusion, condensation phenomena occur \cite{GRVc, GRVk,CCG}, and
there is an analogue of Liggett's comparison inequality which leads to ``propagation of positive correlations'' for appropriately chosen
initial product measures \cite{GRV}.

The fundamental property which makes this model tractable, despite the interaction, is self-duality.
Self-duality, combined with a good coupling of exclusion walkers and independent random walkers has lead in the context of the exclusion process to the proof of the hydrodynamic limit, in the sense of propagation of local equilibrium; see e.g. \cite{MP}, \cite{MIPP} for two good reference on this approach based on $v$-functions.
The construction of a good coupling between exclusion walkers and independent random walkers relies quite strongly on particular properties of
the exclusion process, in particular monotonicity, which in turn leads to a decreasing number
of discrepancies (second class particles) in the courese of time. The SIP is not monotone because of the attractive interaction, and therefore presents new challenges both for coupling and for a strong hydrodynamic limit behavior. By self-duality it is easy to see that the (linear) diffusion equation is
the correct macroscopic equation. However, because the stationary measures are products of discrete Gamma distributions, their tails are exponential and therefore the standard and powerful general method for hydrodynamic limits of gradient systems, namely the GPV (Guo-Papanicolaou-Varadhan)
approach (see \cite{KL}) cannot be applied, neither can we use later developed methods relying on monotonicity. This problem will manifest itself even more  when we pass to the hydrodynamics of the asymmetric inclusion process, where an equation of Burger's type is expected after hyperbolic rescaling of space and time.

In this paper we show that SIP particles and independent random walkers can be coupled such that their distance at time $t$ is $o(\sqrt{t})$ (i.e., goes to zero when divided by $\sqrt{t}$ in the limit $t\to\infty$), which is enough to prove both propagation of local equilibrium and the local equilibrium property of the
boundary driven non-equilibrium steady state. The idea of the coupling is very similar to what happens for simple symmetric exclusion: random walk jumps are performed together and inclusion jumps by inclusion particles only. The discrepancies created on the time scale of potential creation of discrepancies, i.e., when at least two inclusion particles are at nearest neighbor positions, are behaving as a simple random walk. Since inclusion particles in a state
of potential creation of discrepancies will leave this state after an exponential time, with parameter bounded from below uniformly over all
such configurations, we have that the time at which discrepancies can be built up is at most $t^{1/2 +\epsilon}$ with probability
close to one as $t\to\infty$. In this time window the differences between the SIP walkers and corresponding independent random walkers can be decomposed into a sum of quantities behaving as  continuous-time
simple random walk. As a consequence, since the time window is at most $t^{1/2 +\epsilon}$, the total
discrepancy built up will not exceed  $t^{1/4 +\epsilon/2}$ with probability close to one.
This intuition would be correct if there are no effects of multiple particles meeting, i.e., the
symmetry of the created discrepancies are correct for ``two-particle collisions''. We therefore still prove that the effect of higher order collisions, where multiple (more than 2) particles are at neighboring 
positions can be neglected on the time scale we are interested in. 

As a consequence of this coupling result, on the hydrodynamic space and time scale, we can treat
the SIP particles as independent random walkers
and we obtain propagation of local equilibrium with the linear diffusion equation as macro equation.
Furthermore, by applying the coupling to particles absorbed at the boundaries, we obtain the local
equilibrium property of a boundary driven non-equilibrium state obtained by coupling a system
of SIP particles to boundary reservoirs at different densities. Notice that for this non-equilibrium state
no exact solution of the matrix-ansatz type is available.

Our paper is organized
as follows. In sections 2.1--2.2 we introduce the SIP and basic properties such as self-duality and invariant product measures.
In sections 2.3--2.5 we introduce macroscopic profiles, propagation of local equilibrium, the boundary driven SIP, and
the local equilibrium property of its unique non-equilibrium stationary measure. In these sections we use
the coupling between independent walkers and SIP walkers. In section 3 we prove the essential estimate controlling
the distance between the independent walkers and corresponding SIP walkers in this coupling.

\section{Definition and basic properties of the $SIP(m)$}
The simple symmetric inclusion process on $\Z$ with parameter $m$ is defined to be the
Markov process with state space $\Omega=\N^\Z$ (where $\N$ denotes natural numbers including zero) with generator
\be\label{gensip}
\loc^{\rm SIP(m)} f(\eta)= \sum_{i\in \Z} \sum_{j\in \{i-1,i+1\}}\eta_i\left(\frac{m}{2}+\eta_{j}\right) (f(\eta^{i j})-f(\eta))
\ee
where $\eta^{ij}$ denotes the configuration obtained from $\eta\in \Omega$ by removing one particle at $i$
and putting it at $j$: $\eta^{ij}=\eta-\delta_i+\delta_j$ where $\delta_x$ denotes
the configuration with a single particle at location $x\in\Z$ and no particles elsewhere.

The interpretation of the process with generator \eqref{gensip} is that every particle makes the jumps of a continuous-time symmetric nearest neighbor random walk at rate $m/2$  and on top of that there is interaction by inclusion. This interaction is described as follows: to
every pair of particles at nearest neighbor positions is associated a Poisson clock (independent form the others) which rings at rate
one. If this ``invitation clock'' rings for a pair of particles at positions $i,i+1$, then with probability $1/2$ the particle at $i$ joins the particle at $i+1$ or vice versa.

The semigroup associated to the generator \eqref{gensip} is denoted by $S_t$. It is well-defined on multivariate
polynomials
in the variables $\eta_i, i\in \Z$, depending on a finite number of variables. This is a consequence of self-duality, see section \ref{kwakk}
below.
In fact, $S_t f$ is well-defined for many more functions $f$ but we will
need only those (polynomials) here. We denote by $S_t f$ the semigroup working on a function $f$, and $\mu S_t$ for the semigroup working on a probability measure $\mu$ on $\Omega$, i.e., the distribution of $SIP(m)$
at time $t>0$ when started initially from $\mu$.

\subsection{Equilibrium and local equilibrium measures}
The stationary and reversible  product measures
for  the $SIP(m)$ can be easily obtained from a detailed balance computation (cf.\  \cite{GRV}).
If $\nu^{\otimes \Z}$ is a product measure
with marginal $\nu$ then
detailed balance implies
\[
\nu(n)\nu(k) n \left(\frac{m}{2} + k\right)= \nu(n-1)\nu(k+1) (k+1) \left(\frac{m}{2} + n-1\right)
\]
which leads to
\[
\frac{\nu(n)}{\nu(n-1)} \frac{n}{\left(\frac{m}{2} + n-1\right)}= \frac{\nu(k+1)}{\nu(k)}\frac{k+1}{\left(\frac{m}{2} + k\right)}
\]
This gives the discrete Gamma (negative binomial)  distribution $\nu_\lambda$  on $\mathbb{N}$, with scale parameter $\lambda\in[0,1]$
and shape parameter $m/2$,  as possible marginals, i.e. 
\[
\nu(n)=\nu_\lambda (n)= \frac1Z\frac{\lambda^n}{n!}\frac{\Gamma\left(\mee+ n\right)}{\Gamma\left(\mee\right)}, \quad n\in\N.
\]
where $0\leq \lambda<1$ and where
\[
Z= (1-\lambda)^{-m/2}
\]
is the normalizing constant.

We can then introduce ``local'' equilibrium measures, with marginals described above, but
where the constant $\lambda$ is allowed to depend on the location.
For $\lambda :\Z\to [0, 1)$,  we denote by $\nu_\lambda$ the
product measure  on $\Omega$ with probability mass function
\be\label{prod}
\nu_\lambda^{\otimes \Z} (\eta) = \prod_{i\in\Z} \nu_{\lambda(i)}(\eta_i), \quad \eta\in\Omega
\ee
where $\nu_{\lambda(i)}$ is  the discrete Gamma (negative binomial) distribution  on $\N$, with scale parameter $\lambda(i)$ and
shape parameter $m/2$,
defined via the mass function
\be\label{nulaa}
\nu_{\lambda(i)}(n)= \frac{1}{Z_{\lambda(i)}}\frac{\lambda(i)^n}{n!} \frac{\Gamma\left(\mee+ n\right)}{\Gamma\left(\mee\right)}, \quad n\in\N.
\ee
Here $Z_{\lambda(i)}= (1-\lambda(i))^{-m/2}$ is the normalizing constant.
To the ``scale parameter profile'' $\lambda:\Z\to [0,1)$ corresponds the ``density profile''
\be\label{denspro}
\rho(i)= \int \eta_i \nu_\lambda ^{\otimes \Z} (d\eta)= \mee\frac{\lambda(i)}{1-\lambda(i)} .
\ee
For a constant $\lambda$, i.e., $\lambda(i)=\lambda$ for all $i\in\Z$, we of course recover the homogeneous
reversible product measure $\nu_\lambda^{\otimes \Z} (d\eta)$. 
It is natural to expect, but at present not yet proved that these measures are the only ergodic
measures of the SIP.
\subsection{Self-duality}\label{kwakk}
In order to formulate self-duality of the $SIP(m)$, we introduce the self-duality polynomials.
Define
for $k,n \in \N$, $k\leq n$
\be\label{sdpo}
d(k,n)= \frac{n!}{(n-k)!} \frac{\Gamma\left(\mee\right)}{\Gamma\left(\mee+k\right)}
\ee
and for $k>n$, $d(k,n)=0$.
The relation between this polynomial (in $n$) and the discrete Gamma distribution \eqref{nulaa} is given by
\be\label{comodo}
\sum_{n=0}^\infty d(k,n) \nu_{\lambda} (n)=\left( \frac{\lambda}{1-\lambda}\right)^k
\ee

We denote by $\Omega_f$ the set of configurations in $\Omega$ with a finite number of
particles, i.e. $\Omega_f= \{\xi\in \Omega:\sum_{x\in \Z} \xi_x<\infty\}$.
In particular, for $x_1,\ldots, x_n\in \Z$ we denote by
$\sum_{i=1}^n\delta_{x_i}$ the configuration with particles located at positions
$x_1,\ldots, x_n$.

For $\xi\in\Omega_f$
we define
\be\label{sdpoprod}
\caD (\xi, \eta) = \prod_{i\in\Z} d(\xi_i, \eta_i).
\ee
This is well defined since $d(0,n)=1$, for all $n\in\N$.

Self-duality of the $SIP(m)$
is then the following result, proved in \cite{GKRV}, and see also \cite{GRV} for more details
on the self-duality of $SIP(m)$.
\bt
For $\xi\in \Omega_f$, $\eta\in\Omega$ we have
\be\label{selfdualsip}
\E^{\rm SIP(m)}_\eta \caD(\xi, \eta(t)) = \E^{\rm SIP(m)}_\xi \caD (\xi(t), \eta).
\ee
\et
The relation between the product  measures $\nu_\lambda^{\otimes\mathbb{Z}}$ from \eqref{prod} and the duality functions is given by the following
formula \cite{GRV} which easily follows from \eqref{comodo}:
\be\label{measudual}
\int \caD\left(\sum_{i=1}^n\delta_{x_i},\eta\right) \nu_\lambda^{\otimes \Z} (d\eta) = \prod_{i=1}^n \frac{\lambda(x_i)}{1-\lambda(x_i)}.
\ee
In particular, for the homogeneous product measures we have
\be\label{homudal}
\int \caD\left(\sum_{i=1}^n\delta_{x_i},\eta\right) \nu_\lambda^{\otimes \Z} (d\eta) = \left(\frac{\lambda}{1-\lambda}\right)^n.
\ee
For a finite configuration $\xi\in\Omega_f$, we can write
$\xi=\sum_{i=1}^n \delta_{x_i}$ (notice that for different $i$, the corresponding $x_i$ are allowed to be equal), and we will then denote expectation
in the $SIP(m)$ also by $\E_{x_1,\ldots,x_n}^{\rm SIP(m)}$, and the locations
of the corresponding $n$ $SIP(m)$ particles by $X_1(t), \ldots, X_n(t)$.

An important role will be played by the ``correlation functions''
\beq\label{vee}
&&\caV (\nu_\lambda,x_1,\ldots,x_n ; t)
\nonumber\\
&=&\int\E^{\rm SIP(m)}_\eta \caD\left(\sum_{i=1}^n\delta_{x_i},\eta(t)\right) \nu_\lambda^{\otimes \Z}  (d\eta) -
\prod_{i=1}^n \int\E^{\rm SIP(m)}_\eta \caD(\delta_{x_i}, \eta(t)) \nu_\lambda ^{\otimes \Z} (d\eta).\qquad
\eeq
By self-duality, \eqref{measudual} and the fact that a single inclusion walker
is a random walk jumping at rate $\mee$, this can be rewritten as
\beq\label{veedual}
&&\caV (\nu_\lambda^{\otimes \Z} ,x_1,\ldots,x_n ; t)
\nonumber\\
&= &
\E_{x_1,\ldots,x_n}^{\rm SIP(m)} \left(\prod_{i=1}^n \frac{\lambda(X_i(t))}{1-\lambda(X_i(t))}\right)
-
\E_{x_1,\ldots,x_n}^{\rm IRW(m)}\left( \prod_{i=1}^n \frac{\lambda(X_i(t))}{1-\lambda(X_i(t))}\right).
\eeq
Here $IRW(m)$ denote $n$ independent random walkers moving at rate $\mee$, and
$\E_{x_1,\ldots,x_n}^{\rm IRW(m)}$ is the  expectation in this process, initially started at $(x_1,\ldots,x_n)$.
So we see that the function $\caV (\nu_\lambda^{\otimes \Z} ,x_1,\ldots,x_n ; t)$ measures the difference between
expectation over $n$ SIP particles and $n$ random walkers.

Moreover, for a function $\phi:\Z\to\R$, the expectation
\beq
\psi(t,x):= \E_x^{\rm IRW(m)} \phi(X_1(t))
\eeq
solves the discrete
diffusion equation with diffusion constant $\mee$
\be\label{disdif}
\frac{\partial\psi(t,x)}{\partial t} = \frac{m}{2} (\psi(t,x+1)+ \psi(t,x-1)-2\psi(t,x))
\ee
with initial condition $\psi(0,x)=\phi(x)$,
which after suitable scaling  converges to the continuous diffusion equation.

From \eqref{veedual}, as we will see later, it follows that if we can show that in some sense $n$ independent random walkers and $n$ inclusion
walkers can be coupled so that their distance is not too large (as a function of time) then we have, in the sense of
propagation of local equilibrium, that
the inclusion process has the diffusion equation with diffusion constant $m$ \eqref{disdif} as its hydrodynamic limit.

\subsection{Macro profiles and propagation of local equilibrium}
We define a macro profile to be a smooth function
$\pi:\R\to [0,1)$. 
\bd
We say that a sequence $\lambda_N:\Z\to [0,1), N\in\N$ of profiles corresponds
to the macro profile $\pi: \R\to [0,1)$ if for all $y_1,\ldots, y_n\in \R$
\be\label{macro}
\lim_{N\to\infty}\int  \caD\left(\sum_{i=1}^n\delta_{\lfloor Ny_i\rfloor},\eta\right)\nu_{\lambda_N} ^{\otimes \Z} (d\eta)= \prod_{i=1}^n \frac{\pi(y_i)}{1-\pi(y_i)}.
\ee
We denote this property by ``$\lambda_N\approx \pi$''.
\ed
Intuitively this means that around the macro point $y\in \R$, and corresponding micro point $\lfloor Ny\rfloor$, the profile $\lambda_N(\lfloor Ny\rfloor)$
is close to $\pi(y)$, and as  a consequence the product measure $\nu_{\lambda_N }^{\otimes \Z} $ is close
to the product measure $\nu_{\pi}^{\otimes \Z} $.

A simple example of a possible choice for the profiles $\lambda_N$ is given by $\lambda_N(i)= \pi(i/N)$.

\bd
We say that a sequence of probability measures $\mu_N, N\in\N$ on the configuration
space $\Omega$ satisfies the local equilibrium
property (LEP) with profile $\pi$ if for all $y_1,\ldots, y_n\in \R$
\be\label{loceq}
\lim_{N\to\infty}\int  \caD\left(\sum_{i=1}^n\delta_{\lfloor Ny_i\rfloor},\eta\right)\mu_N(d\eta)= \prod_{i=1}^n \frac{\pi(y_i)}{1-\pi(y_i)}.
\ee
We denote this property by $\mu_N\approx {\rm LEP}(\pi)$.
\ed

Because $\nu_{\lambda_N}^{\otimes \Z} $ is a product measure, we have in particular that if a profile $\lambda_N$ corresponds to the macro profile $\pi$, then
$\mu_N=\nu_{\lambda_N}^{\otimes \Z} $ satisfies LEP with macro profile $\pi$, i.e., $\lambda_N\approx \pi$ implies $\nu_{\lambda_N}^{\otimes \Z} \approx {\rm LEP}(\pi)$.

If $\lambda_N\approx\pi$ then we consider the evolution of the local equilibrium
state $\nu_{\lambda_N}^{\otimes \Z} $ after a macroscopic time $N^2t$. By duality, \eqref{vee} and the scaling of simple
random walk to Brownian motion, we have, for $y_1,\ldots, y_n\in \R$
\beq\label{ brom}
&&\int  \E_\eta\caD\left(\sum_{i=1}^n\delta_{\lfloor Ny_i\rfloor},\eta({N^2t})\right)\nu_{\lambda_N}^{\otimes \Z} (d\eta)
\nonumber\\
&=&
\E^{\rm SIP(m)}_{\lfloor Ny_1\rfloor, \ldots, \lfloor Ny_n\rfloor} \left(\prod_{i=1}^n \frac{\lambda_N (X_i(N^2t))}{1-\lambda_N(X_i (N^2t)}\right)
\nonumber\\
&=&\caV( \nu_{\lambda_N}^{\otimes \Z} , \lfloor Ny_1\rfloor, \ldots, \lfloor Ny_n\rfloor, N^2t)
+
\prod_i^n\E^{\rm IRW(m)}_{\lfloor Ny_i\rfloor} \left(\frac{\lambda_N (X_i(N^2t))}{1-\lambda_N(X_i (N^2t)}\right)
\nonumber\\
&=&
\caV( \nu_{\lambda_N}^{\otimes \Z} , \lfloor Ny_1\rfloor, \ldots, \lfloor Ny_n\rfloor, N^2t)
+
\prod_{i=1}^n \psi(t, y_i)
+ o(1)
\eeq
where $o(1)$ goes to zero when $N\to\infty$, and where $\psi(t,y_i)$ is the solution of the diffusion equation
\be\label{difeq}
\frac{\partial\psi(t,y)}{\partial t} =\mee \frac{\partial^2\psi(t,y)}{\partial y^2}
\ee
with initial condition
\be\label{incon}
\psi(0,y)= \frac{\pi(y)}{1-\pi(y)}.
\ee

Therefore, in order to obtain that $\nu_{\lambda_N}^{\otimes \Z} $ after time $N^2t$ satisfies LEP with
macro profile $\psi(t,y)$, it is crucial to prove that the $\caV$ functions vanish in that limit.
More precisely, if we can prove that for all $y_1,\ldots, y_n\in\R$
\be\label{macrovee}
\lim_{N\to\infty} \caV( \nu_{\lambda_N}^{\otimes \Z} , \lfloor Ny_1\rfloor, \ldots, \lfloor Ny_n\rfloor, N^2t)=0
\ee
then we have that $\lambda_N\approx \pi$ implies $\nu_{\lambda_N} ^{\otimes \Z} S_{N^2 t} \approx {\rm LEP} (\psi(t,\cdot))$
where $\psi(t,y)$ satisfies \eqref{difeq}, \eqref{incon}.
We say then that on the macro scale ``\emph{propagation of local equilibrium}'' holds, with corresponding macro equation
\eqref{difeq}.

In turn, because we speed up time by a factor $N^2$, using \eqref{veedual}, \eqref{macrovee} holds if we can find a coupling
$(X_1(t), \ldots, X_n(t);\tilde X_1(t), \ldots,\tilde  X_n(t))$
between $n$ inclusion walkers and $n$ independent random walkers
such that, for all $i$, $|X_i(t)-\tilde X_i(t)|/\sqrt{t}$ converges to zero
in probability as $t\to\infty$.

\subsection{Boundary driven SIP(m)}
The boundary driven $SIP(m)$ on $N\in\mathbb{N}$ lattice sites is a process on
the configuration space $\Omega_N= \N^{\{1,\ldots,N\}}$ defined via its generator
\begin{eqnarray}\label{INC}
{\cal L}^{\rm SIP}f(\eta)&=&{\cal L}_a^{\rm SIP}f(\eta)+{\cal L}_0^{\rm SIP}f(\eta)+{\cal L}_b^{\rm SIP}f(\eta)\\
&=&\alpha (\frac{m}{2}+\eta_1)\left(f(\eta^{0,1})-f(\eta)\right)+\gamma \eta_1 \left(f(\eta^{1,0})-f(\eta)\right)\nonumber\\
&+& \sum_{i=1}^{N-1}\eta_i\left(\mee+\eta_{i+1}\right)\left(f(\eta^{i,i+1})-f(\eta)\right)+ \eta_{i+1}\left(\mee+\eta_i\right)\left(f(\eta^{i+1,i})-f(\eta)\right)\nonumber\\
&+&\sigma (\frac{m}{2}+\eta_N)\left(f(\eta^{N+1,N})-f(\eta)\right)+\beta \eta_N \left(f(\eta^{N,N+1})-f(\eta)\right).\nonumber
\end{eqnarray}
Here, we denote, for $\eta\in\Omega_N$, with slight abuse of notation
\[
\eta^{0,1}= \eta+\delta_1,\  \eta^{1,0}= \eta-\delta_1
\]
and similarly
\[
\eta^{N+1,N}= \eta+\delta_N,\  \eta^{N,N+1}=\eta-\delta_N
\]
this notation turns out to be useful in the dual process, where
we have two more sites, namely $0$ and $N+1$ associated to the reservoirs. Here $\alpha,\beta,\gamma$ and $\sigma$ are positive constants with $\gamma>\alpha$ and $\beta>\sigma$.
In the presence of boundary generators, we have the following duality
result.
\bt
The process $(\eta(t))_{t\ge0}$ defined by \eqref{INC} is dual to the absorbing boundaries process
$(\xi(t))_{t\ge0}$ with configuration space $\Omega_{Dual}=\mathbb{N}_0^{0, 1,\ldots, N, N+1}$ with generator
\begin{eqnarray}
\label{dual-inc}
{\caL}^{\rm SIP}_{\text{Dual}}f(\xi)&=&(\gamma-\alpha) \xi_1 \left(f(\xi^{1,0})-f(\xi)\right)\\
&+& \sum_{i=1}^{N-1}\xi_i(\frac{m}{2}+\xi_{i+1})\left(f(\xi^{i,i+1})-f(\xi)\right)+  \xi_{i+1}(\frac{m}{2}+\xi_i)\left(f(\xi^{i+1,i})-f(\xi)\right)\nonumber\\
&+&(\beta-\sigma) \xi_N \left(f(\xi^{N,N+1})-f(\xi)\right),\;\nonumber
\end{eqnarray}
with duality function
\begin{eqnarray}
\label{dualSIP}
\caD^{\rm SIP}(\xi,\eta)=\left(\frac \alpha{\gamma-\alpha}\right)^{\xi_0} \; \left(\prod_{i=1}^N \frac{\eta_i !}{(\eta_i-\xi_i)!} \, \frac{\Gamma(\frac{m}{2})}{\Gamma(\frac{m}{2}+\xi_i)}\right)\; \left(\frac \sigma {\beta-\sigma}\right)^{\xi_{N+1}}.\;
\end{eqnarray}
\et

\br
Notice that the duality between the operators \eqref{INC} and \eqref{dual-inc}
is valid for all non-negative values of $\alpha,\beta, \gamma, \sigma$.
The restriction $\gamma>\alpha$ and $\beta>\sigma$ ensures that the operator
\eqref{dual-inc} is a Markov generator (i.e., all transition rates are non-negative).
\er

In the dual process the reservoirs have been replaced by absorbing boundaries. This is a considerable simplification, because in the dual process eventually all particles will end up at $0$ or $N+1$.
Therefore, the unique stationary distribution $\nu_{\rho_L, \rho_R}^N $ of the boundary driven process with generator
\eqref{INC} can be completely characterized by the
absorption probabilities of $SIP$ particles absorbed at the boundaries.
This measure $\nu_{\rho_L, \rho_R}^N $
is called ``\emph{non-equilibrium steady state}'' and the non-equilibrium is, in this
case,  created by the boundary reservoirs, i.e., when $\rho_L\not=\rho_R$.

More precisely, 
denoting
$\rho_L= \frac \alpha{\gamma-\alpha},\; \rho_R= \frac\sigma{\beta-\sigma}$ we obtain from
the duality \eqref{dualSIP} the following result for the correlation functions of  $\nu^N_{\rho_L,\rho_R}$ of the boundary driven process with
generator \eqref{INC}.
\be\label{absnes}
\int \caD^{\rm SIP}(\xi,\eta) \nu_{\rho_L, \rho_R}^N (d\eta)= \sum_{k+l=\|\xi\|} \rho_L^k \rho_R^l \pee_{\xi}^{\rm SIP, abs, N} (\xi(\infty)=k\delta_0+l\delta_{N+1})
\ee
where $\|\xi\|=\sum_{x=0}^{N+1} \xi_x$ denotes the total number of dual particles, and where
$\pee^{\rm SIP, abs, N}_\xi $ denotes the path space measure of the dual absorbing SIP process
with generator \eqref{dual-inc} starting from $\xi$.
That is, the correlation functions in the steady state with boundaries are completely determined by the absorption probabilities of
the dual ${\rm SIP}(m)$.

From now on we will choose the parameters
$\gamma= (\rho_L+1) m/2, \alpha= \rho_L m/2, \sigma= \rho_R m/2, \beta= (\rho_R+1)m/2$, which implies that
a single dual particle is a continuous-time random walk  moving at rate $m/2$ and absorbed at $0$ or $N+1$.
Since the absorption probability of such a random walker is linear as a function of the starting point, we obtain as a particular case of \eqref{absnes} the linear
density profile in the steady state with boundaries:
\be\label{lindens}
\int \caD^{\rm SIP}(\delta_i, \eta) \nu_{\rho_L, \rho_R}^{\otimes \Z}  (d\eta)= \rho_L + \frac{(\rho_R-\rho_L) i}{N+1}=: \rho^{R,L}_N(i).
\ee
When we consider this profile on the macro scale we obtain, for $x\in [0,1]$ the linear macro profile
\be\label{lindensmacro}
\lim_{N\to\infty}\int \caD^{\rm SIP}(\delta_{\lfloor xN\rfloor}, \eta) \nu_{\rho_L, \rho_R}^N (d\eta)= \rho_L + (\rho_R-\rho_L)x=: \rho^{R,L, macro}(x).
\ee

\subsection{Local equilibrium property of  steady states with\\ boundaries}
Denote by $\E^{\rm SIP, abs, N}_{x_1,\ldots,x_n}$ the expectation in the process with $n$ SIP particles, moving
on $\{0,\ldots, N+1\}$ according to the generator
\eqref{dual-inc}, and $\E^{\rm IRW, abs, N}_{x_1,\ldots,x_n}$ the corresponding expectation
for independent random walkers. The particles
under $\E^{\rm SIP, abs, N}_{x_1,\ldots,x_n}$ (resp. $\E^{\rm IRW, abs, N}_{x_1,\ldots,x_n}$) move as SIP particles (resp. independent random walkers)
on $\{1,\ldots, N\}$ and are absorbed at the left end $0$ (resp. right end  $N+1$) at rate $\gamma-\alpha$
(resp. $\beta-\sigma$). Absorbed and non-absorbed particles do not interact.
We can then rewrite \eqref{absnes} in this notation:
\beq\label{prodfac}
\int \caD^{\rm SIP}\left(\sum_{i=1}^n \delta_{x_i},\eta\right) \nu_{\rho_L, \rho_R}^N (d\eta)
&=& \sum_{k+l=\|\xi\|} \rho_L^k \rho_R^l \pee^{\rm SIP, abs, N}_{x_1,\ldots, x_n}
(|\{i: X_i(\infty)=0\}|=k)\nonumber\\
&=&\E^{\rm SIP, abs, N}_{x_1,\ldots, x_n}\left(\prod_{i=1}^n \rho_N(X_i(\infty))\right),
\eeq
where we defined $\rho_N(0)= \rho_L, \rho_N(N+1)=\rho_R$.

We can now define the local equilibrium property of $\nu_{\rho_L, \rho_R}^N$ as follows.
Intuitively, this property means that around the macro point $x\in [0,1]$ (corresponding to micro point $\lfloor xN\rfloor$), the measure
$\nu_{\rho_L, \rho_R}^N$ looks like the SIP(m) equilibrium product measure $\nu_{\rho(x)}$ where $\rho(x)$ is the macroscopic linear profile
$\rho^{L,R, macro}$ defined in
\eqref{lindensmacro}.
\bd
Let $\rho: [0,1]\to [0,\infty)$, and $\mu_N$ a collection of probability measures on $\N_0^{\{1,\ldots N\}}$.
The sequence $\mu_N, N\in \N$, satisfies the local equilibrium property with profile $\rho$
(notation $\mu_N \approx{\rm  LEP}(\rho)$ if for all $n\in \N$ and
for all $x_1,\ldots, x_n\in [0,1]$:
\be\label{noneqled}
\lim_{N\to\infty}\int \caD^{\rm SIP}\left(\sum_{i=1}^n \delta_{\lfloor x_iN \rfloor},\eta\right)\mu_N(d\eta)= \prod_{i=1}^n \rho(x_i).
\ee
Here $\N$ is the set of positive integers and $\N_0=\N\cup\{0\}$.
\ed
In view of \eqref{prodfac}
we obtain then the following statement.
If for all $x_1,\ldots,x_n \in [0,1]$ and $a_1, \ldots, a_n\in \{0,N+1\}$
\beq\label{prodabs}
&&\lim_{N\to\infty}\pee^{\rm SIP, abs, N}_{\lfloor x_1N\rfloor,\ldots, \lfloor x_nN\rfloor} (X_1(\infty)=a_1,\ldots, X_n(\infty)= a_n)
\nonumber\\
&=&
\lim_{N\to\infty}\prod_{i=1}^n
\pee^{\rm SIP, abs, N}_{\lfloor x_iN\rfloor} (X_i(\infty)=a_i)
\nonumber\\
&=&
\lim_{N\to\infty}\pee^{I\rm RW, abs, N}_{\lfloor x_1N\rfloor,\ldots, \lfloor x_nN\rfloor} (X_1(\infty)=a_1,\ldots X_n(\infty)= a_n)
\eeq
then $\nu_{\rho_L, \rho_R}^N\approx{\rm LEP}(\rho)$, where
$\rho$ is the macroscopic linear profile
$\rho^{\rm L,R, macro}$ defined in
\eqref{lindensmacro}.
\br
In fact, \eqref{prodabs} is slightly stronger because in 
\eqref{noneqled} the order of particles is not important.
We will in the coupling that we construct always choose to
attach the same label to the random walk particle and
the corresponding $SIP$ particle, and these labels will not
change in the course of the coupling.
\er
In words, asymptotic (as $N\to\infty$) factorization of the absorption probabilities
implies the local equilibrium property of the non-equilibrium steady state.

This factorization \eqref{prodabs} in turn is implied by the existence of a coupling $(X_1(t),\ldots, X_n(t)), (\tilde X_1(t), \ldots, \tilde X_n(t))$
of $n$ SIP(m) particles (with absorption at $0,N+1$)
with $n$ independent random walkers
(with absorption at $0,N+1$) both starting from initial state $(\lfloor x_1N\rfloor, \ldots, \lfloor x_nN\rfloor)$
such that
\be\label{suc}
\lim_{N\to\infty}
\pee \left(X_i(\infty)\not= \tilde X_i(\infty)\ \text{for some}\ i\in \{1,\ldots, n\}\right)=0.
\ee
If we can couple $n$ SIP particles and $n$ independent random walkers on $\Z$
(in a Markovian way) such that
they are at distance $o(\sqrt{t})$ apart when $t\to\infty$, then we can also use that coupling (see Section \ref{S3})  to
achieve \eqref{suc}. Indeed, if the $i^{\rm th}$ SIP particle  is absorbed at $0$, then this happens before time
$bN^2$ for $b$ large enough, with probability close to one. As long as the random walker and the corresponding SIP are not absorbed,
they both move as if they were on the infinite lattice. That is, if the random walk particle is not yet absorbed at the moment that the
SIP particle is absorbed, it is at distance at most $\delta N$ from $0$, with probability close to one, for any $\delta>0$.
Therefore, the probability that it will be absorbed at the other end $N+1$ is at most $\delta$. Hence the probability
that they will be absorbed at different ends tends to zero as $N\to\infty$.

\section{Coupling finitely many SIP particles and \\independent random walkers }\label{S3}
In this section we prove Theorem \ref{th3.1}.
The first non-trivial case to be tackled is the case of
two particles. Then we show that on the time scale $\sqrt{t}$
only ``binary collisions'' play a role, and
the $n$ particle case can thus be reduced to a two particle case.  Results are stated and proved for
the one dimensional case but  easily extend to higher dimenssions using the discussion in Remark
3.1 below.

\subsection{Coupling two SIP particles with two independent random walkers}
We first consider basic coupling $(X(t), Y(t), U(t), V(t))$ such that
$(X(t), Y(t))$ evolves as two $SIP(m)$ particles, $U(t),V(t)$ as two independent
random walkers moving at rate $m/2$.

To define the generator of the basic coupling, we need some notation.
We denote by $e_{13}$, resp.\ $e_{24}$ the vector $(1,0,1,0)$, resp. $(0,1,0,1)$.
Furthermore we denote $\x= (x,y,u,v)$
The generator of the basic coupling then reads
\beq\label{basiccoup2}
&&\loc f(\x)= \mee\sum_{\epsi=\pm 1} \left((f(\x + \epsi e_{13})- f(\x)) + (f(\x + \epsi e_{24})- f(\x))\right)
\nonumber\\
&+&
I(|x-y|=1) (f(x,x,u,v)+f(y,y,u,v)-2 f(x,y,u,v)).
\eeq
This generator is interpreted as follows: random walk jumps are performed together, and
inclusion jumps are performed only in the first two coordinates.
\br
We couple similarly two independent random walk particles absorbed
at left end $0$ and right end $N+1$ and two
$SIP$ particles absorbed
at left end $0$ and right end $N+1$.
The random walk jumps until absorption are performed together, and
the inclusion jumps are performed only by the inclusion particles.
So it can happen that one (e.g. $SIP$)
particle is absorbed and the corresponding
random walk particle not, then the corresponding
random walk particle jumps alone.
\er

The aim is to show that, when the particles are initially at the same location,
i.e., when $X(0)=U(0), Y(0)=V(0)$, then, for all $a>0$, we can keep them closer than $a\,\sqrt{t}$
for $t$ large, with probability close to one. More precisely, we show the following.
\bt\label{th3}
In the basic coupling $(X(t), Y(t), U(t), V(t))$ \eqref{basiccoup2} between two $SIP(m)$ and
two
independent random walkers moving at rate $\mee$, starting at the same initial positions,
we have
\be\label{distance}
\lim_{t\to\infty} \frac{|X(t)-U(t)|^2}{t} =0
\ee
where the limit is in $L^1$, and hence in probability, for every starting position
with $X(0)=U(0), Y(0)=V(0)$. The same statement holds for $|Y(t)-V(t)|$.
\et
\bpr
We are interested in the function $\phi(\x)= x-u$.
An easy computation using \eqref{basiccoup2} gives
\[
\loc \phi = I(y= x+1) - I(y=x-1)
\]
and
\[
\loc (\phi^2) - 2\phi\loc(\phi) = I(y= x+1) + I(y=x-1).
\]
Denoting $z= y-x$ we thus find
\be\label{marti}
\phi (\x(t))-\phi(\x(0)) -\int_0^t I(|z(s)|=1) z(s)= M(t)
\ee
where $\{M(t):t\geq 0\}$ is a martingale with quadratic variation
\be\label{quadva}
<M,M>_t =\int_0^t \left(\loc (\phi^2)(\x(s)) - 2\phi(\x(s))\loc(\phi(\x(s)))\right) ds=\int_0^t I(|z(s)|=1) ds.
\ee
The process $\{z(t):t\geq 0\}$ is a Markov process with generator
\be\label{zgen}
L f(z) = I(|z|=1 ) (f(0)-f(z)) +m\left( f(z+1)+ f(z-1)-2f(z)\right).
\ee
This means   a continuous-time random walk with rate $1+m$ to jump from $1,-1$ towards the origin and
rate $m$ for all other nearest neighbor jumps. This random walk is clearly null recurrent, hence
we conclude from \eqref{quadva} that
\[
\lim_{t\rightarrow\infty}\frac{<M,M>_t}{t}=0.
\]
from which in turn we conclude
\[
\lim_{t\rightarrow\infty}\frac{M(t)}{\sqrt{t}}=0.
\]
Therefore, in order to obtain the desired property \eqref{distance} of the coupled
process $\{\x(t): t\geq 0\}$, we still have to tackle the additive functional
\be\label{addfun}
A(t):=\int_0^t I(|z(s)|=1) z(s) ds.
\ee
We claim now that it is sufficient to show that $\E_0(A(t)^2)/t$ converges to zero when $t\to\infty$, where
$\E_0$ denotes expectation in the process
$\{z(t):t\geq 0\}$ starting at $z(0)=0$. That is, we show that starting at $z(0)=0$ represents no lack of generality.
Indeed, when started at $z(0)=z>0$ to the right of the origin, as long as $0$ is not hit, the random walk
behaves as an ordinary continuous-time nearest neighbor random walk, and therefore, the expectation of
the contribution
to $A(t)$
before hitting zero is dominated by
\[
\left(\int_0^t \E_z(I(|z(s)|=1)ds\right) I(t\leq \tau_0))
\]
where $\tau_0$ denotes the hitting time of $0$. Because $z(s)$, as long as $0$ is not
hit behaves as an ordinary random walk, we have the bound
\[
\int_0^t \E_z(I(|z(s)|=1) ds\leq C \sqrt{t}.
\]
Because $\tau_0$ is finite with probability one, this contribution to $A(t)$ can thus be neglected
(on the relevant time scale $\sqrt{t}$).
Hence, by the strong Markov property, we can then restrict to the case where the starting point
$z(0)=0$.

Let us abbreviate
$p_s(x,y)= \pee (z(s)=y|z(0)=x)$ the transition probability in the Markov process $\{z(s), s\geq 0\}$.
By reflection symmetry and the Markov property, we have
\beq\label{expr}
&&\frac14\E_0(A(T)^2)= \E_0 \left(\int_0^T  I(z(s)=1)\ ds-\int_0^T I(z(s)=-1)\ ds\right)^2
\nonumber\\
&=&
\int_0^T\int_t^T p_t(0,1) p_{s-t} (1,1)\ ds dt-\int_0^T\int_t^T p_t(0,1) p_{s-t} (-1,1)\ dsdt
\nonumber\\
&=&
\int_0^T\int_0^{T-t} p_t(0,1) p_r(1,1) \left(1-\int_0^{T-t-r} f_u(-1,1) du\right) dr\ dt
\eeq
where $f_u(-1,1)$ denotes the probability density for the first hitting time
of $1$, starting from $-1$.

Let us now abbreviate
\[
\psi(T,t,r)=\left(1-\int_0^{T-t-r} f_u(-1,1) du\right)
\]
and for $\delta>0$
\[
\kappa (\delta,T)= \left(1-\int_0^{\delta T} f_u(-1,1) du\right).
\]
By recurrence of the process $\{z(t):t\geq 0\}$,
for all $\delta>0$, $\kappa (\delta, T)\to 0$ as $T\to\infty$.
In fact, even choosing $\delta=\delta_T$ we have
$\lim_{T\to\infty} \kappa (\delta_T, T)\to 0$ as long
as $\delta_T T\to\infty$ as $T\to \infty$ (e.g.
we can choose $\delta_T= T^{-1/2}$).
Because the random walk $z(t)$ moves as an ordinary random
walk outside the set $S=\{-1,0,1\}$, and upon each visit to $S$,  $S$ is  left
after an exponential waiting time, it is straightforward to obtain the bounds
\be\label{bounde}
p_t(0,1)< C/\sqrt{t}, \   p_{t} (1,1)\leq C/\sqrt{t},
\ee
for some constant $C>0$.
These bounds \eqref{bounde}
yield the estimate
\[
\frac{1}{T}\int_0^a\int_0^{b-t} p_t(0,1) p_r(1,1)drdt < C_1ab/T
\]
for some constant $C_1$.

Therefore,
\beq\label{pruts}
&&\frac1T\int_0^T\int_0^{T-t} p_t(0,1) p_r(1,1) \psi(T,r,t) \ drdt
\nonumber\\
&=&
\frac1T\int_0^T\int_0^{T-t} p_t(0,1) p_r(1,1) I(T-t-r> \delta T)\psi(T,r,t)drdt
\nonumber\\
&+&
\frac1T\int_0^T\int_0^{T-t} p_t(0,1) p_r(1,1) I(T-t-r\leq \delta T)\psi(T,r,t)drdt
\nonumber\\
&\leq &
C_2 \kappa (\delta, T)+ C_3\delta.
\eeq
Choosing $\delta=\delta_T= T^{-1/2}$, and taking the limit $T\to\infty$ gives
\be
\lim_{T\to\infty} \frac1T\E_0(A(T)^2)=0
\ee
as we needed.
\epr

\subsection{The general case: Coupling $n$ SIP particles and $n$ independent random walkers }
The aim of this section is to implement the above coupling for $n$ SIP and $n$ random
walkers and show the coupling property
announced in Theorem \ref{th3.1}. To this end we introduce the following notation.  We start by pairing each SIP particle
with one random walker and assigning to each pair a label from the set $E_n=\{1,2,\cdots,n\}$.  Denote by  $\Delta_n=\Z^{E_n}$ the  space of possible positions for the
$n$-particles, i.e.  for $y=(y_1,\cdots,y_n)\in\Delta_n$, $y_i\in\Z$, corresponds to the
position of the $i$th particle. Now for a given particles position configuration  $y\in\Delta_n$,  the number of particles $\eta(x)$ at   site $x\in\Z$,  associted with $y$, is given by $\eta(x)=\sum_{i=1}^n I(y_i=x)$, i.e.  $\eta(x)$ is the  frequency with which $x$
occurs along the sequence $y_1,\ldots,y_n$.

Therefore the coupling  of  interest is that with generator $\mathcal{ L}$ that acts on functions  $f:\Z^n\times\Z^n\mapsto \R$ as
\begin{equation}\label{coupnsipnirw}
\begin{split}
\L f(y,\tilde y)&=\sum_{e=\pm 1}\sum_{i=1}^n \frac m2\left[f(y^{i,i+e},\tilde{y}^{i,i+e})-f(y,\tilde y)\right]\cr
&+\sum_{e=\pm 1}\sum_{i=1}^n\sum_{k=1}^n I(y_k=y_i+e) \left[f(y^{i,i+e},\tilde{y})-f(y,\tilde y)\right],
\end{split}
\end{equation}
where
\begin{equation}
y^{i,i+e}_j=\left\{\begin{array}{ll}
y_j &\mbox{if }\, j\neq i\\
y_i+e &\mbox{otherwise}
\end{array}
\right.
\end{equation}
i.e.  in    $y^{i,i+e}$ each particle, apart from particle $i$, maintains its position, but the
position of the $i$th particle is translated by $e$. Our next result is the general version
of Theorem \ref{th3}.

\begin{theorem}\label{th3.1}
In the coupling $(Y(t),\tilde Y(t))$ \eqref{coupnsipnirw} between the $n$-SIP(m) particles and $n$-independent
random walkers moving at rate $\frac m2$, starting at the same initial positions (i.e. $Y(0)=\tilde Y(0)$), we have
\begin{equation}
\lim_{t\rightarrow\infty }\frac {\left|Y_i(t)-\tilde Y_i(t)\right|^2}{t}=0,\quad \text{for all }\quad1\leq i\leq n,
\end{equation}
where the above limit is in probability.
\end{theorem}

\bpr
As was done for the case with two particles, we consider functions $\phi_i:\Delta_n\times\Delta_n\mapsto \R$  of the form
\begin{equation}
\phi_i(y,\tilde y)=y_i-\tilde y_i, \quad i\in E_n.
\end{equation}
It follows from here that
\begin{equation}\label{Lphi}
\L \phi_i(y,\tilde y)=\sum_{e=\pm1} e\sum_{k=1}^n I(y_k=y_i+e)=
\eta(y_i+1)-\eta(y_i-1),
\end{equation}
where once again  $\eta(x)$ is the number of particles located at the site $x$.  Observe from the second equality of \eqref{Lphi} that $\L\phi_i$ is the difference between the number of
particles located on the right and the left hand side of the position  of the $i$th particle.

Further, we have that
\begin{equation}\label{Lphi1}
\begin{split}
\L (\phi_i^2)(y,\tilde y)-2(\phi\L\phi)(y,\tilde y)&=\sum_{e=\pm1} \sum_{k=1}^n I(y_k=y_i+e)\cr &=\eta(y_i+1)+\eta(y_i-1).
\end{split}
\end{equation}
For each pair $i,k\in E_n$ put $z_{i,k}=y_k-y_i$, i.e., the difference between the positions of
the $i$th and the $k$th SIP(m) particles. This give rise to the  Martingale $M_i$ with
\begin{equation}
M_i(t)=\phi_i(Y(t),\tilde Y(t))-\phi_i(Y(0),\tilde Y(0))- \sum_{k=1}^n\int_0^t I(|z_{i,k}(t)|=1)z_{i,k}(s)ds.
\end{equation}
Here we use the first equality of \eqref{Lphi}. It follows from the first
equality of \eqref{Lphi1} that  the quadratic variation of $M_i$ takes the form
\begin{equation}
\begin{split}
<M_i,M_i>_t&=\int_0^t\left[\L(\phi_i^2)(Y(s),\tilde Y(s))-2(\phi_i\L\phi_i)(Y(s),\tilde Y(s))\right]ds\cr &=\sum_{k=1}^n\int_0^t I(|z_{i,k}(s)|=1)ds=\small{o}(t).
\end{split}
\end{equation}
The  last equality  follows because the process $(Y_1(t),...,Y_n(t))$ is clearly not positive recurrent.

So as before we still have to tackle the analogue of the two particle additive functional
\begin{equation}
A_{i,k}(t)=\int_0^t I(|z_{i,k}(s)|=1)z_{i,k}(s)ds.
\end{equation}
Notice that this resembles very much the two particle quantity: the difference is
however that there are other particles around that could perturb the behavior
of $z_{i,k}(t)$, in particular this is no longer a Markov process now.
We will show however that we can neglect on the relevant time scale the effect
of other particles.
For this we introduce some notation.
Let us call
\begin{equation}
\Delta= \{ s: \exists i\not= k: |z_{i,k}(s)|=1\}
\end{equation}
 the set of \emph{ possible  collision times}:
these are the times at which discrepancies between independent walker
and SIP  particles can arise.
We say that there is a binary collision at time $s\in \Delta$ if
the pair $i,k\in E_n$~ with  $i\not= k$, ~ $|z_{i,k}(s)|=1$~ and $\min\{|z_{i,l}(s)|,|z_{l,k}(s)|\}\geq 2$ for 
all $l\in E_n\setminus\{i,k\}$, i.e., if at time
$s$ precisely two SIP particles are at neighboring positions, and
all other particles are further apart. We call $\caB$ the
set of binary collision times.
We notice that since there are $n$ (i.e., finitely many) particles, the set $\Delta$, when
entered, is left after a random time $T$
which can be dominated by an exponential time $\tau$ with parameter uniformly
bounded from below by $\lambda>0$
(uniform in the realization of the particle configuration).
This is because once entered the set $\Delta$ only a finite number of transitions
is needed to leave it again, and this number is uniformly bounded from
above by an $n$-dependent constant (see figure 1 below for an example of such an exit path).
\begin{figure}[htbp]
\vspace{-.42cm}
\centering
\begin{flushleft}
\includegraphics[scale = 0.56]{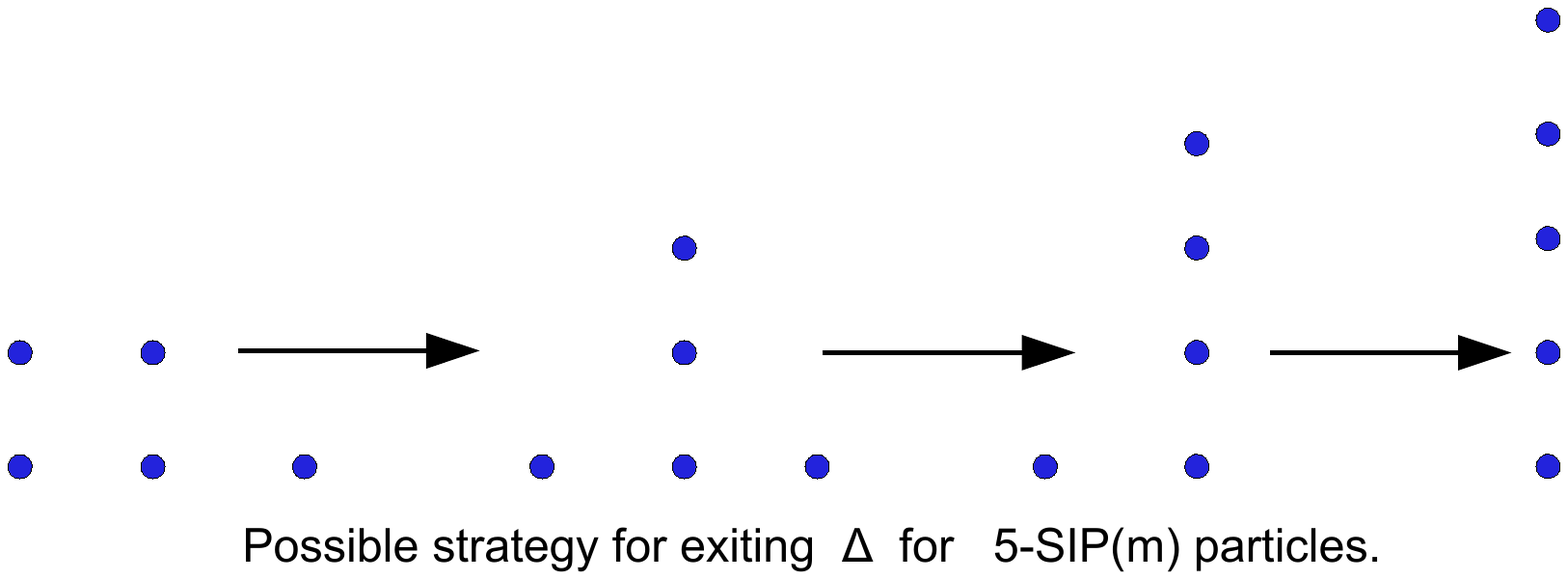}
\end{flushleft}
\vspace{-12.7cm}
\caption{Strategy for exiting $\Delta$}
\label{escaperoute}
\end{figure}
As a consequence, the total time spent in the set $\Delta$ during the time interval $[0,T]$
for the  $n$ SIP particles
has the same
asymptotic leading order behavior as
for $n$ independent random walkers.
Another consequence is that the total time spent
in higher order (than binary) collisions is negligible
on time scale $\sqrt{T}$:
\[
\lim_{T\to\infty} \frac1{\sqrt T}\E_{y}\int_0^T I(s\in \Delta\setminus \caB)ds=0
\]
where the expectation is w.r.t. $n$ SIP particles initially at $y\in \Z^n$.
Indeed, the time that at least three independent random walkers are
at nearest neighbor positions is dominated by $C+ C_n \int_1^T (C/\sqrt{t})^2\leq C+ C_n \ln (T)$.
By the above argument, the same holds  for $n$ SIP particles.
Therefore, if we want to show
\begin{equation}\label{bobobo}
\lim_{t\to\infty}\frac{1}{\sqrt{t}}\E_y A_{i,k}(t)=\lim_{t\to\infty}\frac1{\sqrt{t}}\E_y\int_0^t I(|z_{i,k}(s)|=1)z_{i,k}(s)ds=0
\end{equation}
we can assume that when the $i$ and $k$-th SIP(m) particles create a discrepancy with the
corresponding random walkers because they are at neighboring positions and
make an inclusion jump, the other particles do not interact with these two.
Therefore, we obtain \eqref{bobobo} as in the two particle case.
\epr
\newpage
\br
In this final remarks we indicate some generalizations and further perspectives.
\ben
\item
To generalize to more general finite-range symmetric and translation invariant random walks in $d\geq 1$:
the essential ingredients of the proof of Theorem \ref{th3.1} is that for
the two-particle case discrepancies are created on a time scale of order $\sqrt{t}$ (when
the inclusion walkers are at neighboring positions) and create a difference
of $o(\sqrt{ t})$. Next, multiple collisions can be neglected.
Both statements remain true for general symmetric and translation invariant random walk with bounded
jumps size in $d= 1$.

In $d\geq 2$,
the set $\Delta$ when possible discrepancies occur is then the set of times  when at least two SIP particles  are at distance less than or equal to the maximal jump size.
The fact that in the time interval $[0,t]$, independent
random walkers visit this set $\Delta$ where possible discrepancies are created only during a time window of size $t^{\frac12 +\epsilon}$ with probability close to one as $t\to\infty$ remains identically true.
In fact, in $d=2$ the set $\Delta$ is visited during a time window of order $\log(t)$ and
in $d\geq 3$ of order $1$, so in fact the estimates can even be made stronger, but certainly the result of
Theorem \ref{th3.1} remains valid.
The consequences of the coupling result for the SIP are based on self-duality, and the invariance principle for
a single random walk. For self-duality to hold
it is essential that the underlying random walk is symmetric.
\item
The SIP is dual to an interacting diffusion process with state space $[0,\infty)^\Z$
(energies associated to lattices sites), called the Brownian energy process,
and for integer values of $m$ also to an interacting diffusion process called the Brownian momentum process with state space
$\R^\Z$ (momenta associated to lattice sites), see \cite{GKRV}, \cite{GRV}, and \cite{CGGR} for an overview of all these models.
As a consequence, the results on propagation of local equilibrium and
of the local equilibrium property of the  steady state with boundaries immediately apply to these models.
\een
\er

\end{document}